\documentclass[12pt,oneside]{amsart}
\usepackage{geometry}
\usepackage{parskip}
\usepackage[utf8]{inputenc}
\usepackage{graphicx}
\usepackage{amssymb}
\usepackage{amsmath}
\usepackage{amsthm}
\usepackage{amsfonts}
\usepackage{amsbsy}
\usepackage{hyperref}
\usepackage{xcolor}

\theoremstyle{plain}% default
\newtheorem{thm}{Theorem}

\newtheorem{pop}[thm]{Proposition}
\newtheorem{cor}[thm]{Corollary}

\theoremstyle{definition}
\newtheorem*{defn}{Definition}

\newtheorem*{prf}{Proof}

\theoremstyle{remark}
\newtheorem*{rem}{Remark}

\newcommand{\z}{\mathbb{Z}}

\newcommand{\n}{\mathbb{N}}
\newcommand{\cx}{\mathbb{C}}
\newcommand{\re}{\mathbb{R}}

\author{P. Bhowmik}
\address{Department of Mathematics, University of Rochester, Rochester, NY}
\email{pablo.bhowmik@gmail.com}
\author{S. Deodhar}
\address{Department of Mathematics, University of Rochester, Rochester, NY}
\email{sdeodhar@ur.rochester.edu}
\author{A. Iosevich}
\address{Department of Mathematics, University of Rochester, Rochester, NY}
\email{iosevich@gmail.com}

\thanks{A.~I. was supported in part by the National Science Foundation under NSF DMS-2154232.}

\title{On spectral synthesis in ${\mathbb Z}_N^d$}
\begin{document}
\begin{abstract} A classical result due to Agranovsky and Narayanan \cite{AN04} says that if $f \in L^p({\mathbb R}^n)$ for $2 \leq p \leq \frac{2n}{d}$ and the support of its Fourier transform is contained in a $d$-dimensional $C^1$ submanifold of ${\mathbb R}^n$, where $d<n$, then $f$ is identically equal to $0$. In this paper, we investigate an analogous problem for functions $f: {\mathbb Z}_N^d \to {\mathbb C}$. Bourgain's celebrated result on $\Lambda(p)$ sets \cite{Bou89}, a large-deviation estimate of Hayes \cite{Hay05}, and connections with the theory of exact signal recovery \cite{DS89,MS73,IKLM24,IM24} play an important role.
\end{abstract}

\maketitle
\section{Introduction}
In 2004, Agranovsky and Narayanan \cite{AN04} established the following theorem relating the dimension of the support of the Fourier transform to the function's membership in $L^{p}$.
\begin{thm}[\textit{Agranovsky--Narayanan \cite{AN04}}]\label{AN}
  If $f \in L^{p}(\re^{n})$ and $\operatorname{supp} \hat{f}$ is carried by a $C^{1}$ manifold $M$ of dimension $d < n$, then $f \equiv 0$ provided $2 \leq p \leq \frac{2n}{d}$. If $d=0$, then $f \equiv 0$ for all $2 \leq p < \infty$.
\end{thm}

This result can be viewed as a variant of the Fourier uncertainty principle where the support or concentration condition on the function is replaced by a suitable $L^p$ class. This result was later extended by Senthil Raani \cite{SR14} to sets of fractal dimension, where the packing dimension of the Fourier support determines the integrability threshold. The purpose of this paper is to study the analogous problem for functions on finite abelian groups. 
The central problem we address is the following discrete analogue of the Agranovsky--Narayanan theorem.  
Let \(\{f_N\}_{N}\) be a sequence of functions with \(f_N:\mathbb{Z}_N^d\to\mathbb{C}\) and suppose the Fourier transform \(\widehat f_N\) is supported on a set \(S_N\subset\mathbb{Z}_N^d\) whose size satisfies \(|S_N|\approx N^\alpha\) for some \(\alpha\in(0,d)\).  
We wish to determine conditions on the growth of \(\|f_N\|_{L^p}\) that would force $||f_N||_{\infty}$ to converge to zero as \(N\to\infty\). More precisely, if
$$ {\left( \sum_{x \in {\mathbb Z}_N^d} {|f_{N}(x)|}^p \right)}^{\frac{1}{p}} \leq C$$ for $1 \leq p<\frac{2d}{\alpha}$ and for all $N \in \n$, must $||f_N||_{\infty}$ converge to $0$?

This question is resolved in the affirmative using Theorem~\ref{thm:02} and Corollary~\ref{pop:01} below, and the endpoint case ($p=\frac{2d}{\alpha}$) is discussed in detail. An alternative finite spectral synthesis inequality is introduced in Theorem~\ref{thm:03}. 

Section~\ref{section:bourgain} of this paper is dedicated to establishing the sharpness of the finite spectral synthesis inequalities, where random constructions of sets with optimal Fourier decay and Bourgain's classical result on $\Lambda(p)$ sets play the key role. 

In Section~\ref{section:signalrecovery} we pose and solve the signal recovery problem based on the spectral synthesis results of this paper. More precisely, suppose that, for a sequence of signals $f_N: {\mathbb Z}_N^d \to {\mathbb R}$, the frequencies ${\{\widehat{f}_{N}(m)\}}_{m \in S_N}$ are unobserved, where $S_N \subset {\mathbb Z}_N^d$. Suppose that the  $L^{p}$ norms of $\{f_{N}\}_{N \in \n}$ satisfy a suitable growth bound and $|S_N| \approx N^{\alpha}$. Is it true that $f_N$ can be recovered exactly and uniquely? This question is resolved in Theorem~\ref{thm:04} below, and a recovery algorithm is provided. 
\section{Preliminaries}
We begin by introducing the basic notation and definitions required for our study of spectral synthesis in the setting of finite abelian groups; see \cite{Bab89} for additional background. Throughout this paper, we work primarily on the group
\[
\mathbb{Z}_N^d = \underbrace{\mathbb{Z}_N \times \dots \times \mathbb{Z}_N}_{d\text{ times}},
\]
the \(d\)-fold product of the cyclic group of order \(N\). We identify \(\mathbb{Z}_N\) with \(\{0,1,\dots,N-1\}\) endowed with addition modulo \(N\), and correspondingly view \(\mathbb{Z}_N^d\) as the integer lattice
\[
\{0,1,\dots,N-1\}^d
\]
equipped with componentwise addition modulo \(N\).

For a function \(f:\mathbb{Z}_N^d\to\mathbb{C}\), its \emph{Fourier transform} \(\widehat f:\mathbb{Z}_N^d\to\mathbb{C}\) is defined by
\[
\widehat f(m) = N^{-d/2} \sum_{x\in\mathbb{Z}_N^d} f(x) \, e^{-2\pi i \,x\cdot m/N},\qquad m\in\mathbb{Z}_N^d,
\]
where \(x\cdot m = \sum_{j=1}^d x_j m_j\) denotes the standard inner product.  
The inversion formula reads
\[
f(x)=N^{-d/2}\sum_{m\in\mathbb{Z}_N^d}\widehat f(m)\, e^{2\pi i\, x\cdot m/N},\qquad x\in\mathbb{Z}_N^d.
\]

All \(L^p\) norms in this work are taken with respect to the counting measure; in particular,
\[
\|f\|_{L^p(\mathbb{Z}_N^d)} = \Bigl(\sum_{x\in\mathbb{Z}_N^d}|f(x)|^p\Bigr)^{1/p},\qquad 1\le p<\infty,
\]
and
\[
\|f\|_{L^{\infty}(\mathbb{Z}_N^d)} = \max_{x\in\mathbb{Z}_N^d}|f(x)|.
\]
The same convention is used for functions on the Fourier side. For $1 \leq p \leq \infty$, we denote by $p'$ the conjugate exponent, so that $\frac{1}{p}+\frac{1}{p'}=1$. With this normalization, Plancherel's identity takes the form
\[
\sum_{x\in\mathbb{Z}_N^d}|f(x)|^2 = \sum_{m\in\mathbb{Z}_N^d}|\widehat f(m)|^2.
\]
For a subset $S \subset \mathbb{Z}_N^d$, we denote by $\mathbf{1}_S$ its characteristic function:
\[
\mathbf{1}_S(x) = 
\begin{cases}
1 & \text{if } x \in S, \\
0 & \text{otherwise}.
\end{cases}
\]

\vskip.25in 

\section{Integrability conditions and spectral synthesis inequalities in discrete settings}
In this section, we present a quantitative version of the Agranovsky--Narayanan theorem and establish analogous results in the discrete setting of $\mathbb{Z}_N^d$. A subsequent quantitative refinement of their theorem relates the support of the Fourier transform to the growth of the $L^p$ norm over expanding balls. In what follows, $B_r$ denotes the open ball of radius $r$ centered at the origin, and $\mathbf{1}_{B_r}$ is its characteristic function.
\begin{cor}[\textit{Quantitative version of Theorem~\ref{AN}}]
  Let $f \in L^1_{\mathrm{loc}}(\mathbb{R}^n)$ and suppose $\operatorname{supp} \widehat{f}$ is contained in a $C^1$ manifold $M$ of dimension $d<n$. Assume that $2 \leq p < \frac{2n}{d}$ and that there exists $\delta$ with
  \[
  0 \leq \delta < \frac{1}{2}\left(\frac{2n}{p}-d\right)
  \]
  such that
    \[
    \frac{\|f \mathbf{1}_{B_r}\|_p}{r^{\delta}\log r}=o(1)
    \quad \text{as } r \to \infty.
    \]
    Then $f \equiv 0$.
\end{cor}
This result follows directly from the arguments in \cite[Theorem 1]{AN04}, adapted to the growth condition on the $L^{p}$ norm of the function. We refer the reader to \cite{AN04} for the details. 
The above result provides a quantitative relationship between the growth of the $L^{p}$ norm and the support of its Fourier transform. In the discrete setting, one obtains analogous estimates for sequences of functions $f_{N}: \mathbb{Z}_{N}^{d} \rightarrow \mathbb{C}$ whose Fourier transforms are supported on sets of size $N^{\alpha}$, where $\alpha \in (0,d)$. No compatibility between the functions for different values of $N$ is required. We now proceed to state our primary result.
\begin{thm}\label{thm:02}
If $f: \z_{N}^{d} \rightarrow \cx$ is a function such that $\operatorname{supp}\widehat{f} \subset S$, where $S$ is a subset of $\z_{N}^{d}$, and if $2 \leq p < \infty$, then
    \begin{equation}\label{eq:02}
        ||f||_{L^{\infty}({\mathbb Z}_N^d)} \leq \sqrt{\frac{|S|}{N^{\frac{2d}{p}}}} \cdot {||f||}_{L^p({\mathbb Z}_N^d)}.\end{equation} 
\end{thm}
 The following corollary of Theorem~\ref{thm:02} yields an Agranovsky--Narayanan-type integrability condition for the sequence of functions $\{f_{N}\}$ on $\mathbb{Z}_{N}^{d}$. Here we consider the Fourier transform to have support  of size approximately \( N^{\alpha} \), denoted by \( |\operatorname{supp}(\widehat{f})| \approx N^{\alpha} \). This means that \(C_{1} N^{\alpha} \leq |\operatorname{supp}(\widehat{f})| \leq C_{2} N^{\alpha} \), where \( C_{1}, C_{2} \) are constants independent of \( N \).

\begin{cor}\label{pop:01}
     Let $\{f_N\}_{N\in\mathbb{N}}$ be a sequence of functions
$f_N:\mathbb{Z}_N^d\to\mathbb{C}$. Suppose that there exists a constant
$C>0$, independent of $N$, such that
$$
\|f_N\|_{L^p(\mathbb{Z}_N^d)}\leq C
$$
for every $N\in\mathbb{N}$. Let
$
S_N:=\operatorname{supp}(\widehat{f_N}),
$
and assume that $|S_N|\approx N^\alpha$ for some $\alpha\in(0,d)$. If
$
1\leq p<\frac{2d}{\alpha},
$
then
\[
\|f_N\|_{L^\infty(\mathbb{Z}_N^d)}\longrightarrow 0
\qquad\text{as }N\to\infty.
\]
 \end{cor}
 \begin{prf}
  Suppose first that $p \geq 2$. Substituting the uniform bound $C$ into (\ref{eq:02}), we obtain
  $$||f_{N}||_{{L^{\infty}(\z_{N}^{d})}} \leq |S_{N}|^{\frac{1}{2}}N^{-\frac{d}{p}}C.$$
  Since $p < \frac{2d}{\alpha}$, the right-hand side converges to $0$.

  Suppose now that $1 \leq p < 2$. Fourier inversion, the Cauchy--Schwarz inequality, and Plancherel's theorem give
  $$
  ||f_N||_{L^{\infty}(\z_N^d)}
  \leq N^{-\frac{d}{2}}|S_N|^{\frac{1}{2}}||f_N||_{L^2(\z_N^d)}.
  $$
  Since counting measure is being used, $||f_N||_{L^2} \leq ||f_N||_{L^p}$. Consequently,
  $$
  ||f_N||_{L^{\infty}(\z_N^d)}
  \lesssim N^{-\frac{d}{2}}N^{\frac{\alpha}{2}},
  $$
  which converges to $0$ because $\alpha<d$.
  \qed
 \end{prf}
 Theorem~\ref{thm:02} is intrinsic to $\mathbb{Z}_N^d$ and gives a quantitative relationship between the $L^p$ norm of a function and the size of its Fourier support.
 A key difference between the continuous and discrete settings lies in the behavior at the endpoint exponent. In the Euclidean case, the vanishing result holds inclusively at the endpoint $p = \frac{2n}{d}$ (where $n$ is the ambient dimension and $d$ is the dimension of the supporting submanifold), forcing any such function to be identically zero. However, in the discrete setting, this property fails at the analogous endpoint $p = \frac{2d}{\alpha}$; that is, there exists a sequence of functions $\{f_{N}\}$ such that $||f_{N}||_{L^p(\z_{N}^{d})}$ remains uniformly bounded at this exponent, yet $||f_{N}||_{L^{\infty}(\z_{N}^{d})}$ does not converge to zero.
In the next section, we will provide an example of such functions using large deviation estimates. Another distinction between the continuous and discrete settings appears in the range
$1\leq p<2$. Theorem~\ref{thm:02} is stated for $p\geq 2$ and therefore yields no
direct conclusion in the subquadratic range. In the discrete setting, however,
Corollary~\ref{pop:01} continues to give a vanishing result for $1\leq p<2$. This
extension is not a new spectral-synthesis phenomenon; rather, it follows from the
monotonicity of $\ell^p$ norms with respect to counting measure on
$\mathbb{Z}_N^d$, namely,
\[
\|f\|_{L^2(\mathbb{Z}_N^d)}
\leq
\|f\|_{L^p(\mathbb{Z}_N^d)}
\qquad\text{whenever }1\leq p<2.
\]
Combining this observation with the estimate at $p=2$ yields the conclusion in the
subquadratic range. Thus, this additional range depends on the discrete
counting-measure normalization rather than on an intrinsic improvement of the
spectral-synthesis inequality. We now turn to another spectral synthesis-type inequality.
 \begin{thm}\label{thm:03}
     Let $2 \leq p \leq \infty$. If $f: \z_{N}^{d} \rightarrow \cx$ is a function such that $\operatorname{supp}\widehat{f} \subset S$, where $S$ is a subset of $\z_{N}^{d}$, then
     \begin{equation}\label{eq:03}
||f||_{L^{\infty}({\mathbb Z}_N^d)} \leq N^{-\frac{d}{2}} \cdot {||f||}_{L^p({\mathbb Z}_N^d)} \cdot {||\widehat{\mathbf{1}_S}||}_{L^{p'}({\mathbb Z}_N^d)}.
\end{equation}
 \end{thm}
 \begin{rem}
Theorem~\ref{thm:02} can be recovered from Theorem~\ref{thm:03} using Hölder's inequality and Plancherel's theorem as follows:
\begin{align*}
\|\widehat{\mathbf{1}_S}\|_{L^{p'}(\mathbb{Z}_N^d)} &\leq N^{d\left(\frac{1}{p'} - \frac{1}{2}\right)} \|\widehat{\mathbf{1}_S}\|_{L^2(\mathbb{Z}_N^d)} && \text{(Hölder's inequality).} \\
\intertext{Using Plancherel's theorem, we obtain}
\|\widehat{\mathbf{1}_S}\|_{L^2(\mathbb{Z}_N^d)} &= \|\mathbf{1}_S\|_{L^2(\mathbb{Z}_N^d)} = |S|^{\frac{1}{2}}. \\
\intertext{Substituting this into the inequality, we obtain}
\|\widehat{\mathbf{1}_S}\|_{L^{p'}(\mathbb{Z}_N^d)} &\leq N^{d\left(\frac{1}{p'} - \frac{1}{2}\right)} |S|^{\frac{1}{2}}.
\end{align*}
Substituting this into (\ref{eq:03}), we obtain Theorem~\ref{thm:02}.
\end{rem}
The large-scale behavior of Fourier coefficients of characteristic functions is a well-studied phenomenon. For a nonzero function $f$, Theorem~\ref{thm:03} can be rearranged to give a lower bound on the $L^{p'}$ norm of the Fourier transform of a characteristic function. Consequently, identifying functions $f$ for which equality holds in the above bound is of particular interest, and we will investigate these cases in the next section.
  \section{A large deviation principle, Bourgain's Theorem and sharpness of spectral synthesis inequalities} \label{section:bourgain}
In this section, we first introduce a large deviation inequality and use it to construct an extremal example demonstrating the sharpness of Theorem~\ref{thm:02}. We then present a result due to Bourgain and discuss its application in proving the optimality of the spectral synthesis inequality introduced earlier.
\subsection{Random sets and sharpness example}
 In the previous section, we established an integrability result for sequences of functions on $\mathbb{Z}_{N}^{d}$ with increasing $N$. We now prove the sharpness of Corollary~\ref{pop:01}. Specifically, for $p \geq \frac{2d}{\alpha}$, we construct a sequence of functions $\{f_{N}\}$ whose $L^{p}$ norms are uniformly bounded and for which $||f_{N}||_{L^{\infty}(\z_{N}^{d})}$ does not converge to 0. We begin by stating a corollary of an Azuma-type large deviation inequality due to Hayes.
 \begin{defn}
      Let $G$ be a finite abelian group, and let $f : G \to \mathbb{C}$ be any function. We denote
$$
\Phi(f) := \max_{\chi} \left| \sum_{x \in G} \chi(x) f(x) \right|,
$$
where the maximum is taken over all characters $\chi : G \to \mathbb{C}^\times$ \textit{except} the principal character $\chi \equiv 1$. In the special case where $f : G \to \{0,1\}$ is the characteristic function of a subset $S \subseteq G$, we will write $\Phi(S)$ in place of $\Phi(f)$.
\end{defn}
Using an Azuma-type large deviation inequality, Hayes proved the following theorem.
\begin{thm}[\textit{Hayes \cite{Hay05}}]\label{thm:06}
    Let $G$ be a finite abelian group of order $n$, let $0<k<n$, and set $m=\min\{k,n-k\}$. If $S'$ is chosen uniformly at random among all subsets of $G$ of cardinality $k$, then, for every $a>0$,
    \begin{equation}\label{eq:04}
        \mathbb{P}\left(\Phi(S') \geq a\right) \leq 2ne^2\exp\left(-\frac{a^2}{8m}\right).
    \end{equation}
\end{thm}
When $k\leq \frac{n}{2}$, substituting $a=k^{\frac{1}{2}}\log k$ into (\ref{eq:04}) shows that, provided $k$ grows polynomially in $n$, the largest non-trivial Fourier coefficient is bounded by $k^{\frac{1}{2}}\log k$ for almost all subsets of cardinality $k$. Using this result, we will construct a sequence of sets of a given size in $\mathbb{Z}_{N}^{d}$ such that the functions
\[
f_{N} = \frac{N^{d/2}}{|S_{N}|}\widehat{\mathbf{1}_{S_{N}}}
\] 
have uniformly bounded $L^{p}$ norms for  $p > \frac{2d}{\alpha}$, and $||f_{N}||_{L^{\infty}(\z_{N}^{d})}$ does not converge to $0$ as $N \rightarrow \infty$.
\begin{pop}\label{pop:02}
     Let $\alpha \in (0,d)$. If $p >\frac{2d}{\alpha}$, then there exists a sequence of functions $f_{N}: \z_{N}^{d} \rightarrow \cx$ with $|\operatorname{supp}(\widehat{f_{N}})| \approx N^{\alpha}$ such that
     $$||f_{N}||_{L^p(\z_{N}^{d})} \leq C \quad \mathrm{and} \quad ||f_{N}||_{L^{\infty}(\z_{N}^{d})} > C_{0}$$
     for all $N \in \n$, where $C$ and $C_{0}$ are constants independent of $N$.
 \end{pop}

 \begin{prf}
    Define the following function:
    $$f_{N}(x) = \frac{N^{d/2}}{|S_{N}|}\widehat{\mathbf{1}_{S_{N}}}(x)$$
    where $S_N \subset \z_{N}^{d}$ is a set of size $|S_N| \approx N^\alpha$ to be specified later. Next, we define the maximum non-trivial Fourier coefficient of the set:
    $$\phi(S_N) = \max\{N^{d/2}|\widehat{\mathbf{1}_{S_N}}(x)| ~:~ x \neq 0\}.$$ 
    First, we establish a deterministic bound on the $L^p$ norm of $f_N$ by splitting the $L^p$ sum into the zero frequency and the non-zero frequencies. At the origin, evaluating the Fourier transform yields exactly:
$$f_N(0) = \frac{N^{d/2}}{|S_N|} \widehat{1_{S_N}}(0) = \frac{N^{d/2}}{|S_N|} \left( N^{-d/2} |S_N| \right) = 1.$$
For all $x \neq 0$, the definition of $\phi(S_N)$ guarantees that $|f_N(x)| \leq \phi(S_N)|S_N|^{-1}$. Summing over the group elements, we obtain:
\begin{equation}\label{eq:10}
    ||f_N||_{L^p(\mathbb{Z}_{N}^{d})} = \left( |f_N(0)|^p + \sum_{x \neq 0} |f_N(x)|^p \right)^{\frac{1}{p}} \leq \left( 1 + N^d \left( \phi(S_N)|S_N|^{-1} \right)^p \right)^{\frac{1}{p}}.
\end{equation}
    
    To find a suitable set $S_N$, we consider all sets of a fixed size $|S_N|$ and apply Theorem~\ref{thm:06} to obtain the large deviation bound:
    $$\mathbb{P}(\phi(S_N) \leq |S_N|^{1/2+\epsilon}) \geq 1-cN^{d}e^{-|S_N|^{2\epsilon}/8},$$
    where
    $$
    0<\epsilon<\frac{1}{2}-\frac{d}{\alpha p}.
    $$
    
    Since $|S_N| \approx N^\alpha$ and $\alpha<d$, we have $|S_N|<\frac{N^d}{2}$ for all sufficiently large $N$. Moreover, the error term $cN^{d}e^{-|S_N|^{2\epsilon}/8}$ becomes strictly less than $1$. Thus, $\mathbb{P}(\phi(S_N) \leq |S_N|^{1/2+\epsilon}) > 0$, guaranteeing that we can choose a specific set $S_{N} \subset \z_{N}^{d}$ that satisfies this bound.
    
    Substituting this set $S_N$ into \eqref{eq:10}, we get:
    $$||f_N||_{L^p(\z_{N}^{d})} \leq \left(N^d|S_N|^{(\epsilon-\frac{1}{2})p}+1\right)^{\frac{1}{p}}.$$
    Because $\epsilon>0$ is sufficiently small and $p > \frac{2d}{\alpha}$, the exponent of $N$ becomes negative. Therefore, for large enough $N$, we have $N^d|S_N|^{(\epsilon-\frac{1}{2})p} < 1$, which yields:
    $$||f_N||_{L^p(\z_{N}^{d})} \leq 2^{\frac{1}{p}}.$$
    
    To bound the $L^\infty$ norm, it is easy to see that evaluating $f_N$ at the origin gives:
    $$||f_{N}||_{L^{\infty}(\z_{N}^{d})} \geq f_N(0) = 1.$$ 
    
    This establishes the required bounds for all sufficiently large $N$. For smaller values of $N$, any arbitrary choice of $S_{N} \subset \z_{N}^{d}$ yields finite, non-zero norms. By taking the maximum bounds over these finitely many  cases, we obtain the desired universal constants $C$ and $C_0$.
\qed
 \end{prf}
In the Euclidean setting, it is known that if the Fourier transform of a function is supported on a submanifold of codimension greater than one, then the minimal exponent $p$ for which a non-zero function $f \in L^p(\mathbb{R}^n)$ exists is typically strictly greater than $\frac{2n}{d}$, where $d$ denotes the dimension of the supporting submanifold. A notable example is the moment curve \cite{GIZZ23}, for which the sharp integrability threshold is known to be $\frac{n^2 + n + 2}{2}$, substantially larger than $2n$ in general.

On the other hand, in the discrete setting where $\alpha$ can take any value in $(0,d)$, the integrability result is sharp for sets of all sizes. We also observe that the above proposition yields a sharpness example only for the case $p > \frac{2d}{\alpha}$. In what follows, we present a theorem due to Bourgain and construct a sharpness example for the endpoint case $p = \frac{2d}{\alpha}$.
\subsection{Bourgain's Theorem and applications}
Let us present Bourgain's theorem from his $\Lambda(p)$ paper \cite{Bou89}, which will be essential for constructing the sharpness example. In what follows, $\|\cdot\|_{p}$ will denote the $L^{p}$ norm with respect to the standard uniform probability measure on the given space.
\begin{thm}[\textit{Bourgain \cite{Bou89}}]
    Let $\{\phi_i\}_{i\in\z_N^d}$ be the family of characters on the group $\z_{N}^{d}$, and let $2<p< \infty$. Then there is a subset $S$ of $\z_{N}^{d}$ with $|S| \approx N^{\frac{2d}{p}}$ such that
    \begin{equation}\label{eq:15}
       ||\sum_{i\in S} a_{i}\phi_{i}||_{p} \leq C(p)\left(\sum_{i \in S} |a_{i}|^{2}\right)^{\frac{1}{2}}
    \end{equation}
    for all scalars $\{a_{i}\}$, where $C(p)$ depends only on $p$.
\end{thm}
The preceding theorem is a special case of a fundamental result proved by Bourgain for orthogonal systems of functions on locally compact groups.

\begin{pop}\label{pop:03}
    If $p = \frac{2d}{\alpha}$, then there exists a sequence of functions $f_{N}: \z_{N}^{d} \rightarrow \cx$ such that $|\operatorname{supp}(\widehat{f_{N}})| \approx N^{\alpha}$ and
     $$||f_{N}||_{L^p(\z_{N}^{d})} \leq C ~~\mathrm{and}~~ ||f_{N}||_{{L^{\infty}(\z_{N}^{d})}} >C_{0}$$
    for all $N \in \n$, where $C$ and $C_{0}$ are constants independent of $N$.
\end{pop}
\begin{prf}
    Here we follow the same template as the proof of Proposition~\ref{pop:02}. Let us first define the function $f_{N}$ as
    $$f_{N}(x) = \frac{N^{d/2}}{|S_{N}|}\widehat{\mathbf{1}_{S_{N}}}(x),$$
    where we choose $S_{N} \subset \z_{N}^{d}$ to be a set of size $|S_{N}| \approx N^{\frac{2d}{p}}$ satisfying \eqref{eq:15}. For such a set $S_{N}$, we have
    $$||\widehat{\mathbf{1}_{S_{N}}}||_{L^p(\z_{N}^{d})} \leq C(p) N^{\frac{d}{p}}N^{-\frac{d}{2}}|S_{N}|^{\frac{1}{2}}.$$
    Since $p = \frac{2d}{\alpha}$, we know that $|S_N|^{\frac{1}{2}} \approx N^{\frac{d}{p}}$. Substituting this yields
    $$||\widehat{\mathbf{1}_{S_{N}}}||_{L^p(\z_{N}^{d})} \lesssim C(p) |S_{N}|N^{-\frac{d}{2}}.$$
    Recalling the definition of $f_N$, this implies $||f_{N}||_{L^p(\z_{N}^{d})} \leq C(p)$. Furthermore, evaluating $f_N$ at the origin gives
    $$||f_{N}||_{L^{\infty}(\z_{N}^{d})} \geq 1.$$
    This completes the proof.
    \qed
\end{prf}
The proof of the preceding proposition demonstrates that Bourgain's theorem allows us to select a set of size $\sim N^{\frac{2d}{p}}$ with optimal $L^p$ control. We will apply this key idea to construct a sharpness example for inequality~(\ref{eq:03}).
  \begin{pop}\label{pop04}
    Let $\alpha \in (0, d)$ and let $p$ satisfy $2 \leq p \leq \frac{2d}{\alpha}$.
    For all sufficiently large $N$, there exist a set $S \subset \mathbb{Z}_N^d$ with $|S| \approx N^{\alpha}$ and a function $f: \mathbb{Z}_N^d \to \mathbb{C}$ with $\operatorname{supp} \widehat{f} \subset S$ such that
    \[
    \|f\|_{L^{\infty}(\mathbb{Z}_N^d)} \; \approx \; N^{-d/2} \, \|f\|_{L^{p}(\mathbb{Z}_N^d)} \, \|\widehat{\mathbf{1}_S}\|_{L^{p'}(\mathbb{Z}_N^d)},
    \]
    where the implied constants depend only on $d$ and $p$.
\end{pop}
\begin{prf}
    From (\ref{eq:03}), it suffices to show that there is a set $S$ satisfying the given condition such that
    \begin{equation}\label{eq:12}
    ||f||_{L^p(\z_{N}^{d})}||\widehat{\mathbf{1}_S}||_{L^{p'}(\z_{N}^{d})} \lesssim N^{d/2}||f||_{{L^{\infty}(\z_{N}^{d})}}.
    \end{equation}
    
    We will use $\Lambda(p)$ sets supplied by Bourgain's theorem to establish the sharpness. First, we obtain a deterministic bound on $||\widehat{\mathbf{1}_S}||_{L^{p'}(\z_{N}^{d})}$. Since $p' \leq 2$, by Hölder's inequality we have
    $$||\widehat{\mathbf{1}_S}||_{L^{p'}(\z_{N}^{d})} \leq ||\widehat{\mathbf{1}_S}||_{L^{2}(\z_{N}^{d})}N^{(\frac{d}{p'}-\frac{d}{2})}.$$
    Using Plancherel's theorem, we can conclude that
    $$||\widehat{\mathbf{1}_S}||_{L^{2}(\z_{N}^{d})} = ||\mathbf{1}_S||_{L^{2}(\z_{N}^{d})} = |S|^{\frac{1}{2}},$$
    thus,
    \begin{equation}\label{eq:16}
    ||\widehat{\mathbf{1}_S}||_{L^{p'}(\z_{N}^{d})} \leq |S|^{\frac{1}{2}}N^{(\frac{d}{p'}-\frac{d}{2})}.    
    \end{equation}
    
    We now choose a set $S$ of size approximately $N^{\alpha}$ with the required $L^p$ control. If $p=2$, we may take any such set, since the desired estimate follows immediately from Plancherel's theorem. Suppose that $p>2$. Bourgain's theorem gives a set $S_{1}$ such that
    $$
    |S_{1}| \approx N^{\frac{2d}{p}}
    $$
    and
    $$||\sum_{i\in S_{1}} a_{i}\phi_{i}||_{p} \leq C(p)\left(\sum_{i \in S_{1}} |a_{i}|^{2}\right)^{\frac{1}{2}}$$
    for all scalar sequences $(a_i)_{i\in S_1}$. Since $p \leq \frac{2d}{\alpha}$, we have $N^{\alpha}\leq N^{\frac{2d}{p}}$. We may therefore choose a subset $S\subset S_1$ with $|S|\approx N^{\alpha}$. The $\Lambda(p)$ estimate is inherited by subsets. Taking $a_i=N^{-\frac{d}{2}}$ for $i\in S$ and $a_i=0$ otherwise gives
    $$||\widehat{\mathbf{1}_S}||_{p} \leq C(p) N^{-\frac{d}{2}}|S|^{\frac{1}{2}}.$$
    
    Since we used the normalized $L^p$ norm in Bourgain's theorem, after denormalizing, we obtain
    \begin{equation}\label{eq:17}
    ||\widehat{\mathbf{1}_S}||_{L^p(\z_{N}^{d})} \leq C(p) N^{\frac{d}{p}}N^{-\frac{d}{2}}|S|^{\frac{1}{2}}.     
    \end{equation}
    
    Thus, for this particular set, combining (\ref{eq:16}) and (\ref{eq:17}), we obtain
    $$||\widehat{\mathbf{1}_S}||_{L^p(\z_{N}^{d})}||\widehat{\mathbf{1}_S}||_{L^{p'}(\z_{N}^{d})} \leq C(p)\left(N^{\frac{d}{p}}N^{-\frac{d}{2}}|S|^{\frac{1}{2}}\right) \cdot \left(N^{\frac{d}{p'}}N^{-\frac{d}{2}}|S|^{\frac{1}{2}}\right) = C(p)|S|.$$
    
    Define
    $$
    f(x):=\widehat{\mathbf{1}_S}(-x).
    $$
    Then $f$ is the inverse Fourier transform of $\mathbf{1}_S$, so that $\widehat f=\mathbf{1}_S$. In particular,
    $$\operatorname{supp}(\widehat{f})=S \quad \mathrm{and} \quad ||f||_{{L^{\infty}(\z_{N}^{d})}} = \frac{|S|}{N^{\frac{d}{2}}}.$$
    Substituting the respective values and using the above equation, we can conclude that $f$ satisfies \eqref{eq:12}, giving us the desired sharpness.
    \qed
\end{prf}
In all the above examples, we constructed random sets to illustrate sharpness cases. In the following proposition, we show that coordinate subgroups of \( \mathbb{Z}_{N}^{d} \) also provide examples of functions that attain equality in Theorem~\ref{thm:03}. This helps to partially establish sharpness for \( p > \frac{2d}{\alpha} \).
\begin{pop}\label{pop:05}
For any integers \(d \ge 2\), \(N \ge 1\), and \(k\) with \(0 \le k \le d\), let
\[
H=\mathbb{Z}_N^k\times\{0\}^{d-k}.
\]
Then there exists a function $f: \z_{N}^{d} \rightarrow \mathbb{C}$ whose Fourier transform is supported on the subgroup $H$, which has size \(|H| = N^k\), such that
    \[
    \| f \|_{L^{\infty}(\mathbb{Z}_N^d)} = N^{-d/2} \| f \|_{L^p(\mathbb{Z}_N^d)} \| \widehat{\mathbf{1}_H} \|_{L^{p'}(\mathbb{Z}_N^d)},
    \]
    for every \(2 \le p \le \infty\) (with \(p'\) the conjugate exponent).
\end{pop}
\begin{prf}
    Let $H=\mathbb{Z}_N^k\times\{0\}^{d-k}$ and consider the function $f = \widehat{\mathbf{1}_H}$. Since applying the Fourier transform twice reflects a function, we have
    $$
    \widehat f(m)=\mathbf{1}_H(-m)=\mathbf{1}_H(m).
    $$
    Thus, $\operatorname{supp}(\widehat f)=H$. We also have
     $$\widehat{\mathbf{1}_H}(m) = N^{-d/2}\sum_{x \in H}\chi(x \cdot m),$$
     where $\chi(n) = e^{\frac{-2\pi i n}{N}}$. Let $H^{\perp}$ denote the orthogonal complement of $H$. If $m \in H^{\perp}$, then, by definition, we have $\chi(m \cdot x) = 1$ for all $x \in H$. Thus,
     \begin{equation}\label{eq:H}
     \widehat{\mathbf{1}_H}(m) = N^{-d/2} \sum_{x \in H} 1 = N^{-d/2}|H| = N^{-d/2} N^{k}.
     \end{equation}

     If $m \notin H^\perp$, then the character $x \mapsto \chi(m \cdot x)$ is non-trivial on $H$, and hence its sum over $H$ vanishes:
     $$\widehat{\mathbf{1}_H}(m) = N^{-d/2} \sum_{x \in H} \chi(x \cdot m) = 0.$$
     Hence $\mathrm{supp}(\widehat{\mathbf{1}_H}) = H^{\perp}$. We can infer from this observation and equation \eqref{eq:H} that
     $$||\widehat{\mathbf{1}_H}||_{L^p(\z_{N}^{d})} = |H^{\perp}|^{1/p} N^{-d/2} N^{k},$$
     and
     $$||\widehat{\mathbf{1}_H}||_{L^{p'}(\z_{N}^{d})} = |H^{\perp}|^{1/p'} N^{-d/2} N^{k}.$$
     Combining the two equations, we obtain
     $$||\widehat{\mathbf{1}_H}||_{L^p(\z_{N}^{d})} \, ||\widehat{\mathbf{1}_H}||_{L^{p'}(\z_{N}^{d})} = |H^{\perp}| N^{-d} N^{2k}.$$
     As $|H^{\perp}| = N^{d - k}$, we obtain
     $$||\widehat{\mathbf{1}_H}||_{L^p(\z_{N}^{d})} \, ||\widehat{\mathbf{1}_H}||_{L^{p'}(\z_{N}^{d})} = N^{k}.$$
     Since $\widehat{\mathbf{1}_H}$ is supported on $H^{\perp}$ and has constant value $N^{-d/2} N^{k}$ from \eqref{eq:H}, we conclude
     $$||\widehat{\mathbf{1}_H}||_{{L^{\infty}(\z_{N}^{d})}} = N^{-d/2} N^{k}.$$
     Thus, for $f = \widehat{\mathbf{1}_H}$, we have
     $$N^{-d/2} ||f||_{L^p(\z_{N}^{d})} \, ||\widehat{\mathbf{1}_H}||_{L^{p'}(\z_{N}^{d})} = ||f||_{L^{\infty}(\z_{N}^{d})}.$$
     \qed
\end{prf}
In general, if a set $S$ satisfies
$$|\mathrm{supp}(\widehat{\mathbf{1}_S})| \cdot |S| = N^{d},$$
that is, if it attains equality in the uncertainty principle, then one can arrive at a similar equality as in the preceding proposition. The proof follows along similar lines.
  \section{Spectral synthesis and Signal Recovery} \label{section:signalrecovery}
    In this section, we shift our focus to the problem of signal recovery. In their seminal work, Donoho and Stark \cite{DS89} framed the uncertainty principle in terms of the support sizes of a function and its Fourier transform, subsequently deriving signal recovery conditions based on these support sizes. In our study so far, we have established spectral synthesis inequalities and presented an uncertainty principle in terms of the $L^p$ norm and the size of the Fourier support. Accordingly, we will proceed to study the signal recovery problem with constraints in terms of $L^p$ growth and the size of the support on the Fourier side. 
     
     Given a sequence of signals $\{f_N\}_{N\in\mathbb{N}}$ whose $L^p$ norms
satisfy a suitable growth bound, under what conditions can we recover
$f_N$ exactly and uniquely for all sufficiently large $N$? We may think of $f_N$ as the information available about the signal at resolution $N$. To make the problem tractable, we impose a structural assumption on the sets $S_N$ of missing frequencies and a quantization assumption on the signals. Specifically, we impose a restriction on the range of admissible values for the signal, as follows. For $\delta>0$, define
\[
D_{\delta,N}:=\left\{u:\mathbb{Z}_N^d\to\mathbb{R}: u(x)\in\delta\mathbb{Z}\text{ for every }x\in\mathbb{Z}_N^d\right\}.
\]
This $\delta$-separation condition can be interpreted as a quantization constraint, implying that the function takes values in the lattice $\delta \mathbb{Z}$. We now formulate a precise recovery theorem.
\begin{thm}\label{thm:04}
Let $\alpha\in(0,d)$, let $2\leq p<\frac{2d}{\alpha}$, and let $\delta>0$. Suppose that $f_{N}\in D_{\delta,N}$ is a sequence of quantized signals and that the Fourier transform of $f_{N}$ is transmitted except for the frequencies $\{\widehat{f}_N(m)\}_{m\in S_N}$, which are unobserved. Suppose that there is a constant $C_{\mathrm{size}}>0$ such that
\[
|S_N|\leq C_{\mathrm{size}}N^{\alpha}.
\]
Then $f_{N}$ is the unique solution of the following minimization problem if
\begin{equation}\label{eq:011}
    ||f_{N}||_{L^{p}(\mathbb{Z}_{N}^{d})}<\frac{\delta}{2\sqrt{C_{\mathrm{size}}}}N^{(\frac{d}{p}- \frac{\alpha}{2})}
\end{equation}
for all sufficiently large $N$:
\begin{equation}\label{eq:11}
    f_{N} = \operatorname*{argmin}_{u\in D_{\delta,N}}||u||_{L^p(\z_{N}^{d})}
    \quad \text{subject to} \quad
    \widehat{u}(m)=\widehat{f}_{N}(m) \text{ for every }m \notin S_N.
    \end{equation}
\end{thm}
\begin{rem}
The assumption $f_N\in D_{\delta,N}$ means that the signal takes values in the known lattice $\delta\mathbb{Z}$, although the value at each point is unknown. The theorem above guarantees the unique recovery of such signals, subject to the assumptions on the transmitted Fourier transform and the $L^p$ growth of $f_N$. The result above readily extends to the case of fine quantization, where the separation depends on $N$, i.e., $\delta_N \sim N^{-\gamma}$ for some $\gamma \ge 0$. In this regime, the corresponding nontrivial recoverability range is
\[
2\leq p<\frac{2d}{\alpha+2\gamma}.
\]
\end{rem}
\section{Proofs of theorems}
\subsection{Proof of Theorem~\ref{thm:02}}
By Fourier inversion and the assumption that $\widehat{f}$ is supported in $S$,  
 $$ f(x)=N^{-\frac{d}{2}} \sum_{m \in S} \chi(x \cdot m) \widehat{f}(m).$$
 It follows that 
 $$ |f(x)| \leq N^{-\frac{d}{2}} \cdot {|S|}^{\frac{1}{2}} {\left( \sum_{m \in {\mathbb Z}_N^d} {|\widehat{f}(m)|}^2 \right)}^{\frac{1}{2}}.$$
 By Plancherel, this quantity is equal to 
 $$ N^{-\frac{d}{2}} \cdot {|S|}^{\frac{1}{2}} {\left( \sum_{x \in {\mathbb Z}_N^d} {|f(x)|}^2 \right)}^{\frac{1}{2}}.$$
% $$={|S|}^{\frac{1}{2}} {\left( N^{-d} \sum_{x \in {\mathbb Z}_N^d} {|f(x)|}^2 \right)}^{\frac{1}{2}}.$$ 
By Hölder's inequality, this quantity is bounded by 
 $$ {|S|}^{\frac{1}{2}} {\left( N^{-d} \sum_{x \in {\mathbb Z}_N^d} {|f(x)|}^p \right)}^{\frac{1}{p}}$$
 $$ =\sqrt{\frac{|S|}{N^{\frac{2d}{p}}}} \cdot {||f||}_{L^p({\mathbb Z}_N^d)}.$$
 \qed
\subsection{Proof of Theorem~\ref{thm:03}}
     To prove (\ref{eq:03}), we write 
$$ \widehat{f}(m)=\widehat{f}(m)1_S(m).$$

\vskip.125in 

Let $\widetilde{g}(x)=g(-x)$ and define convolution by
$$
(u*v)(x)=\sum_{y\in\mathbb{Z}_N^d}u(y)v(x-y).
$$
Since the inverse Fourier transform of $\mathbf{1}_S$ is $\widetilde{\widehat{\mathbf{1}_S}}$, it follows that
$$ f(x)=N^{-\frac{d}{2}} \cdot \left(f*\widetilde{\widehat{\mathbf{1}_S}}\right)(x).$$

\vskip.125in 

By Hölder's inequality, we conclude that 
$$ |f(x)| \leq N^{-\frac{d}{2}} \cdot {||f||}_{L^p({\mathbb Z}_N^d)} \cdot {||\widehat{\mathbf{1}_S}||}_{L^{p'}({\mathbb Z}_N^d)}.$$
\qed
\subsection{Proof of Theorem~\ref{thm:04}}
     Let $g_N$ be any minimizer in (\ref{eq:11}). Such a minimizer exists because $f_N$ is admissible and the $L^p$ norm is coercive on the discrete set $D_{\delta,N}$. Define $h_{N}:=f_{N}-g_{N}$. By Minkowski's inequality,
\[
||h_{N}||_{L^p(\z_{N}^{d})}\leq||f_{N}||_{L^p(\z_{N}^{d})}+||g_{N}||_{L^p}\leq 2||f_{N}||_{L^p(\z_{N}^{d})}.
\]
 The support of $\widehat{h_{N}}$ is contained in $S_{N}$ since $\widehat{f_{N}}$ and $\widehat{g_{N}}$ agree away from $S_{N}$. Since both $f_N$ and $g_N$ belong to $D_{\delta,N}$, if $h_{N}$ is not identically zero, then
\[
||h_{N}||_{{L^{\infty}(\z_{N}^{d})}}\geq\delta.
\]

Applying Theorem~\ref{thm:02} and the observations above, we see that
\[
\delta\leq||h_{N}||_{{L^{\infty}(\z_{N}^{d})}}
\leq 2\sqrt{C_{\mathrm{size}}}N^{\frac{\alpha}{2}-\frac{d}{p}}||f_{N}||_{L^{p}(\mathbb{Z}_{N}^{d})}.
\] 

It follows that, if we assume
\[ 
||f_{N}||_{L^{p}(\mathbb{Z}_{N}^{d})}<\frac{\delta}{2\sqrt{C_{\mathrm{size}}}}N^{(\frac{d}{p}- \frac{\alpha}{2})}
\]
we obtain a contradiction and conclude that $h_{N}$ must be identically $0$. Since $g_N$ was an arbitrary minimizer, $f_N$ is the unique minimizer. \qed
  \begin{rem}
The condition $p < \frac{2d}{\alpha}$ is essential for a nontrivial uniform recovery statement. If $p > \frac{2d}{\alpha}$, the recovery bound on the right-hand side of (\ref{eq:011}) decays to zero as $N \to \infty$. Since any non-zero quantized signal satisfies $||f_N||_{L^p} \geq \delta$, this regime would force $f_N \equiv 0$ for sufficiently large $N$.
\end{rem}

\end{document}